\documentclass[11pt, a4paper]{article}
\usepackage[cp1251]{inputenc}
\usepackage{amsmath} \usepackage{euscript} \usepackage{marvosym}
\usepackage[dvipsnames]{xcolor}
\usepackage{mymatrx6}
\usepackage{graphicx}
\usepackage{longtable}
\usepackage{mathrsfs}
\usepackage{amssymb, amscd}
\usepackage{adigraph}
\Linestrue\Autonumtrue
\usepackage{mathrsfs}
\def\bluecommand{\textcolor[cmyk]{1,0,0,0}}
\def\redcommand{\textcolor[cmyk]{0,1,1,0}}
\def\greencommand{\textcolor[cmyk]{1,0,1,0}}
\def\yellowcommand{\textcolor[cmyk]{0,0.25,1,0}}
\def\graycommand{\textcolor[cmyk]{0,0,0,0.2}}

\begin{document}

\newcounter{bnomer} \newcounter{snomer}
\newcounter{bsnomer}
\setcounter{bnomer}{0}
\renewcommand{\thesnomer}{\thebnomer.\arabic{snomer}}
\renewcommand{\thebsnomer}{\thebnomer.\arabic{bsnomer}}
\renewcommand{\refname}{\begin{center}\large{\textbf{References}}\end{center}}

\setcounter{MaxMatrixCols}{14}

\newcommand\restr[2]{{
  \left.\kern-\nulldelimiterspace 
  #1 
  \right|_{#2} 
}}

\newcommand{\sect}[1]{%
\setcounter{snomer}{0}\setcounter{bsnomer}{0}
\refstepcounter{bnomer}
\par\bigskip\begin{center}\large{\textbf{\arabic{bnomer}. {#1}}}\end{center}}
\newcommand{\sst}[1]{%
\refstepcounter{bsnomer}
\par\bigskip\textbf{\arabic{bnomer}.\arabic{bsnomer}. {#1}}\par}
\newcommand{\defi}[1]{%
\refstepcounter{snomer}
\par\medskip\textbf{Definition \arabic{bnomer}.\arabic{snomer}. }{#1}\par\medskip}
\newcommand{\theo}[2]{%
\refstepcounter{snomer}
\par\textbf{Theorem \arabic{bnomer}.\arabic{snomer}. }{#2} {\emph{#1}}\hspace{\fill}$\square$\par}
\newcommand{\mtheop}[2]{%
\refstepcounter{snomer}
\par\textbf{Theorem \arabic{bnomer}.\arabic{snomer}. }{\emph{#1}}
\par\textsc{Proof}. {#2}\hspace{\fill}$\square$\par}
\newcommand{\mcorop}[2]{%
\refstepcounter{snomer}
\par\textbf{Corollary \arabic{bnomer}.\arabic{snomer}. }{\emph{#1}}
\par\textsc{Proof}. {#2}\hspace{\fill}$\square$\par}
\newcommand{\mtheo}[1]{%
\refstepcounter{snomer}
\par\medskip\textbf{Theorem \arabic{bnomer}.\arabic{snomer}. }{\emph{#1}}\par\medskip}
\newcommand{\theobn}[1]{%
\par\medskip\textbf{Theorem. }{\emph{#1}}\par\medskip}
\newcommand{\theoc}[2]{%
\refstepcounter{snomer}
\par\medskip\textbf{Theorem \arabic{bnomer}.\arabic{snomer}. }{#1} {\emph{#2}}\par\medskip}
\newcommand{\mlemm}[1]{%
\refstepcounter{snomer}
\par\medskip\textbf{Lemma \arabic{bnomer}.\arabic{snomer}. }{\emph{#1}}\par\medskip}
\newcommand{\mprop}[1]{%
\refstepcounter{snomer}
\par\medskip\textbf{Proposition \arabic{bnomer}.\arabic{snomer}. }{\emph{#1}}\par\medskip}
\newcommand{\theobp}[2]{%
\refstepcounter{snomer}
\par\textbf{Theorem \arabic{bnomer}.\arabic{snomer}. }{#2} {\emph{#1}}\par}
\newcommand{\theop}[2]{%
\refstepcounter{snomer}
\par\textbf{Theorem \arabic{bnomer}.\arabic{snomer}. }{\emph{#1}}
\par\textsc{Proof}. {#2}\hspace{\fill}$\square$\par}
\newcommand{\theosp}[2]{%
\refstepcounter{snomer}
\par\textbf{Theorem \arabic{bnomer}.\arabic{snomer}. }{\emph{#1}}
\par\textsc{Sketch of the proof}. {#2}\hspace{\fill}$\square$\par}
\newcommand{\exam}[1]{%
\refstepcounter{snomer}
\par\medskip\textbf{Example \arabic{bnomer}.\arabic{snomer}. }{#1}\par\medskip}
\newcommand{\deno}[1]{%
\refstepcounter{snomer}
\par\textbf{Notation \arabic{bnomer}.\arabic{snomer}. }{#1}\par}
\newcommand{\lemm}[1]{%
\refstepcounter{snomer}
\par\textbf{Lemma \arabic{bnomer}.\arabic{snomer}. }{\emph{#1}}\hspace{\fill}$\square$\par}
\newcommand{\lemmp}[2]{%
\refstepcounter{snomer}
\par\medskip\textbf{Lemma \arabic{bnomer}.\arabic{snomer}. }{\emph{#1}}
\par\textsc{Proof}. {#2}\hspace{\fill}$\square$\par\medskip}
\newcommand{\coro}[1]{%
\refstepcounter{snomer}
\par\textbf{Corollary \arabic{bnomer}.\arabic{snomer}. }{\emph{#1}}\hspace{\fill}$\square$\par}
\newcommand{\mcoro}[1]{%
\refstepcounter{snomer}
\par\textbf{Corollary \arabic{bnomer}.\arabic{snomer}. }{\emph{#1}}\par\medskip}
\newcommand{\corop}[2]{%
\refstepcounter{snomer}
\par\textbf{Corollary \arabic{bnomer}.\arabic{snomer}. }{\emph{#1}}
\par\textsc{Proof}. {#2}\hspace{\fill}$\square$\par}
\newcommand{\nota}[1]{%
\refstepcounter{snomer}
\par\medskip\textbf{Remark \arabic{bnomer}.\arabic{snomer}. }{#1}\par\medskip}
\newcommand{\propp}[2]{%
\refstepcounter{snomer}
\par\medskip\textbf{Proposition \arabic{bnomer}.\arabic{snomer}. }{\emph{#1}}
\par\textsc{Proof}. {#2}\hspace{\fill}$\square$\par\medskip}
\newcommand{\hypo}[1]{%
\refstepcounter{snomer}
\par\medskip\textbf{Conjecture \arabic{bnomer}.\arabic{snomer}. }{\emph{#1}}\par\medskip}
\newcommand{\prop}[1]{%
\refstepcounter{snomer}
\par\textbf{Proposition \arabic{bnomer}.\arabic{snomer}. }{\emph{#1}}\hspace{\fill}$\square$\par}

\newcommand{\proof}[2]{%
\par\medskip\textsc{Proof{#1}}. \hspace{-0.2cm}{#2}\hspace{\fill}$\square$\par\medskip}

\makeatletter
\def\iddots{\mathinner{\mkern1mu\raise\p@
\vbox{\kern7\p@\hbox{.}}\mkern2mu
\raise4\p@\hbox{.}\mkern2mu\raise7\p@\hbox{.}\mkern1mu}}
\makeatother

\newcommand{\okr}[2]{%
\refstepcounter{snomer}
\par\medskip\textbf{{#1} \arabic{bnomer}.\arabic{snomer}. }{\emph{#2}}\par\medskip}

\newcommand{\Ind}[3]{%
\mathrm{Ind}_{#1}^{#2}{#3}}
\newcommand{\Res}[3]{%
\mathrm{Res}_{#1}^{#2}{#3}}
\newcommand{\epsi}{\varepsilon}
\newcommand{\tri}{\triangleleft}
\newcommand{\Supp}[1]{%
\mathrm{Supp}(#1)}
\newcommand{\NSupp}[1]{%
\mathrm{NSupp}(#1)}
\newcommand{\SSu}[1]{%
\mathrm{SingSupp}(#1)}

\newcommand{\lee}{\leqslant}
\newcommand{\gee}{\geqslant}
\newcommand{\reg}{\mathrm{reg}}
\newcommand{\Dyn}{\mathrm{Dyn}}
\newcommand{\Ann}{\mathrm{Ann}\,}
\newcommand{\Cent}[1]{\mathbin\mathrm{Cent}({#1})}
\newcommand{\PCent}[1]{\mathbin\mathrm{PCent}({#1})}
\newcommand{\Irr}[1]{\mathbin\mathrm{Irr}({#1})}
\newcommand{\Exp}[1]{\mathbin\mathrm{Exp}({#1})}
\newcommand{\empr}[2]{[-{#1},{#1}]\times[-{#2},{#2}]}
\newcommand{\sreg}{\mathrm{sreg}}
\newcommand{\ilm}{\varinjlim}
\newcommand{\wdth}{\mathrm{wd}}
\newcommand{\plm}{\varprojlim}
\newcommand{\codim}{\mathrm{codim}\,}
\newcommand{\GKdim}{\mathrm{GKdim}\,}
\newcommand{\chara}{\mathrm{char}\,}
\newcommand{\rk}{\mathrm{rk}\,}
\newcommand{\chr}{\mathrm{ch}\,}
\newcommand{\Ker}{\mathrm{Ker}\,}
\newcommand{\id}{\mathrm{id}}
\newcommand{\Ad}{\mathrm{Ad}}
\newcommand{\Gh}{\mathrm{Gh}}
\newcommand{\col}{\mathrm{col}}
\newcommand{\row}{\mathrm{row}}
\newcommand{\high}{\mathrm{high}}
\newcommand{\low}{\mathrm{low}}
\newcommand{\pho}{\hphantom{\quad}\vphantom{\mid}}
\newcommand{\fho}[1]{\vphantom{\mid}\setbox0\hbox{00}\hbox to \wd0{\hss\ensuremath{#1}\hss}}
\newcommand{\wt}{\widetilde}
\newcommand{\wh}{\widehat}
\newcommand{\ad}[1]{\mathrm{ad}_{#1}}
\newcommand{\tr}{\mathrm{tr}\,}
\newcommand{\GL}{\mathrm{GL}}
\newcommand{\SL}{\mathrm{SL}}
\newcommand{\SO}{\mathrm{SO}}
\newcommand{\Or}{\mathrm{O}}
\newcommand{\Sp}{\mathrm{Sp}}
\newcommand{\SuppD}{\mathbb{S}\mathrm{upp}}
\newcommand{\Sa}{\mathrm{S}}
\newcommand{\Sing}{\mathrm{Sing}}
\newcommand{\Ua}{\mathrm{U}}
\newcommand{\Andre}{\mathrm{Andre}}
\newcommand{\Aord}{\mathrm{Aord}}
\newcommand{\Mat}{\mathrm{Mat}}
\newcommand{\Stab}{\mathrm{Stab}}
\newcommand{\htt}{\mathfrak{h}}
\newcommand{\spt}{\mathfrak{sp}}
\newcommand{\slt}{\mathfrak{sl}}
\newcommand{\sot}{\mathfrak{so}}

\newcommand{\vfi}{\varphi}
\newcommand{\aad}{\mathrm{ad}}
\newcommand{\vpi}{\varpi}
\newcommand{\teta}{\vartheta}
\newcommand{\Bfi}{\Phi}
\newcommand{\Fp}{\mathbb{F}}
\newcommand{\Rp}{\mathbb{R}}
\newcommand{\Zp}{\mathbb{Z}}
\newcommand{\Cp}{\mathbb{C}}
\newcommand{\Ap}{\mathbb{A}}
\newcommand{\Pp}{\mathbb{P}}
\newcommand{\Kp}{\mathbb{K}}
\newcommand{\Np}{\mathbb{N}}
\newcommand{\ut}{\mathfrak{u}}
\newcommand{\at}{\mathfrak{a}}
\newcommand{\glt}{\mathfrak{gl}}
\newcommand{\hei}{\mathfrak{hei}}
\newcommand{\nt}{\mathfrak{n}}
\newcommand{\kt}{\mathfrak{k}}
\newcommand{\mt}{\mathfrak{m}}
\newcommand{\rt}{\mathfrak{r}}
\newcommand{\rad}{\mathfrak{rad}}
\newcommand{\bt}{\mathfrak{b}}
\newcommand{\unt}{\underline{\mathfrak{n}}}
\newcommand{\gt}{\mathfrak{g}}
\newcommand{\vt}{\mathfrak{v}}
\newcommand{\pt}{\mathfrak{p}}
\newcommand{\Xt}{\mathfrak{X}}
\newcommand{\Po}{\mathcal{P}}
\newcommand{\PV}{\mathcal{PV}}
\newcommand{\Uo}{\EuScript{U}}
\newcommand{\Fo}{\EuScript{F}}
\newcommand{\Do}{\EuScript{D}}
\newcommand{\Eo}{\EuScript{E}}
\newcommand{\Jo}{\EuScript{J}}
\newcommand{\Iu}{\mathcal{I}}
\newcommand{\Mo}{\mathcal{M}}
\newcommand{\Nu}{\mathcal{N}}
\newcommand{\Ro}{\mathcal{R}}
\newcommand{\Co}{\mathcal{C}}
\newcommand{\Ko}{\mathcal{K}}
\newcommand{\So}{\mathcal{S}}
\newcommand{\Lo}{\mathcal{L}}
\newcommand{\Ou}{\mathcal{O}}
\newcommand{\Uu}{\mathcal{U}}
\newcommand{\Tu}{\mathcal{T}}
\newcommand{\Au}{\mathcal{A}}
\newcommand{\Vu}{\mathcal{V}}
\newcommand{\Du}{\mathcal{D}}
\newcommand{\Bu}{\mathcal{B}}
\newcommand{\Sy}{\mathcal{Z}}
\newcommand{\Sb}{\mathcal{F}}
\newcommand{\Gr}{\mathcal{G}}
\newcommand{\Xu}{\mathcal{X}}
\newcommand{\Op}{\mathbb{O}}
\newcommand{\chv}{\mathrm{chv}}
\newcommand{\rtc}[1]{C_{#1}^{\mathrm{red}}}

\newcommand{\JSpec}[1]{\mathrm{JSpec}\,{#1}}
\newcommand{\MSpec}[1]{\mathrm{MSpec}\,{#1}}
\newcommand{\PSpec}[1]{\mathrm{PSpec}\,{#1}}
\newcommand{\APbr}[1]{\mathrm{span}\{#1\}}
\newcommand{\APbre}[1]{\langle #1\rangle}
\newcommand{\APro}[1]{\setcounter{AP}{#1}\Roman{AP}}\newcommand{\apro}[1]{{\rm\setcounter{AP}{#1}\roman{AP}}}
\newcommand{\ot}{\xleftarrow[]{}}
\newcounter{AP}


\author{Stepan Bondar\and Mikhail Ignatev}
\date{}
\title{Tangent cones to Schubert varieties for Kac--Moody groups}\maketitle
\begin{abstract} Let $G$ be the affine Kac--Moody group of type $\wt A_{n-1}$, $B$ be an Iwahori subgroup in $G$, $\Fo=G/B$ be the flag variety, and $W$ be the Weyl group of $G$. Given distinct involutions $w_1,~w_2\in W$, we prove that the tangent cones $C_{w_1}$, $C_{w_2}$ to the corresponding Schubert subvarieties $X_{w_1}$ and $X_{w_2}$ of $\Fo$ at the point $p=e\mod B$ do not coincide as subvarieties of the tangent space to $\Fo$ at the point $p$. This generalizes similar results in the finite-dimensional setting. The main technical tool we used are combinatorics of the embeddings of the Weyl groups of different ranks and coadjoint orbis for the unipotent radical of the group $B$.

\end{abstract}


\sect{Introduction and the main result}\label{sect:intro}\addcontentsline{toc}{subsection}{\ref{sect:intro}. Introduction}

\let\thefootnote\relax\footnote{The work was supported by RSF (project no. 25--11--00302).}

Let $G$ be a complex reductive algebraic group, $T$ a maximal torus in~$G$, $B$ a Borel subgroup in~$G$ containing $T$, and $U$ the unipotent radical of $B$. Let $\Phi$ be the root system of $G$ with respect to~$T$, $\Phi^+$ the set of positive roots with respect to $B$, $\Delta$ the set of simple roots, and $W$ the Weyl group of $\Phi$ (see \cite{Bourbaki}, \cite{Humphreys} and \cite{Humpreys2} for basic facts about algebraic groups and root systems).

Denote by $\Fo=G/B$ the flag variety and by $X_w\subseteq\Fo$ the Schubert subvariety corresponding to an element $w$ of the Weyl group $W$. Denote by $\Ou=\Ou_{p,X_w}$ the local ring at the point $p=eB\in X_w$. Let $\mt$ be the maximal ideal of~$\Ou$. The decreasing sequence of ideals $$\Ou\supseteq\mt\supseteq\mt^2\supseteq\ldots$$ is a filtration on $\Ou$. We define $R$ to be the graded algebra $$R=\mathrm{gr}\,\Ou=\bigoplus_{i\geq0}\mt^i/\mt^{i+1}.$$ By definition, the \emph{tangent cone} $C_w$ to the Schubert variety $X_w$ at the point $p$ is the spectrum of~$R$: $C_w=\mathrm{Spec}\,R$. Obviously, $C_w$ is a subvariety of the tangent space $T_pX_w\subseteq T_p\Fo$. Indeed, let $\Au$ be the symmetric algebra of the vector space $\mt/\mt^2$, or, equivalently, the algebra of regular functions on the tangent space $T_pX_w$. Since $R$ is generated as $\Cp$-algebra by $\mt/\mt^2$, it is a~quotient ring $R=\Au/I$. A hard problem in studying geometry of $X_w$ is to describe $C_w$ \cite[Chapter 7]{BilleyLakshmibai}.

In 1962, A. Kirillov invented the orbit method, which allows to classify unitary irreducible representations of the unipotent radical $U$ of the group $B$ in terms of the coadjoint action of $U$ on the dual space $\ut^*$ to the Lie algebra $\nt$ of $U$ (see \cite{Kirillov62} and \cite{Kirillov04}. (Coadjoint action is by definition the dual action to the usual adjoint action of $U$ on $\ut$.)

In 2013, D. Eliseev and A. Panov computed tangent cones $C_w$ for all $w\in W$ in the case\break $G=\mathrm{SL}_n(\mathbb{C})$, $n\leq5$ \cite{EliseevPanov}. Using their computations, A. Panov formulated the following conjecture, which plays an important role in the view  of the Kirillov's orbit method \cite{Kirillov03}.

\hypo{\textup{(Panov, 2013)} Let $w_1$\textup{,} $w_2$ be \label{mconj}involutions\textup{,} i.e.\textup{,} $w_1^2=w_2^2=\id$. If $w_1\neq w_2$\textup{,} then $C_{w_1}\neq C_{w_2}$ as subvarieties of $T_p\Fo$.}

One can easily check that it is enough to prove the conjecture for irreducible root systems. In 2013, D. Eliseev and the second author proved this conjecture in types $A_n$, $F_4$ and $G_2$ \cite{EliseevIgnatyev}. In \cite{BochkarevIgnatyevShevchenko}, M. Bochkarev, A. Shevchenko and the second author proved the conjecture in types $B_n$ and~$C_n$. In \cite{IgnatyevShevchenko1}, A. Shevchenko and the second author proved that the conjecture is true if $\Phi$ is of type $D_n$ and $w_1$, $w_2$ are so-called basic involutions. In the paper \cite{IgnatyevShevchenko20}, A. Shevchenko and the second author proved that the conjecture is true for so-called good pairs of involutions for $\Phi=E_6$, $E_7$ and $E_8$.

Now, let $G$ be a Kac--Moody group (see, e.g., \cite{Kumar02} for precise definitions). One can define $B$, $T$, $\Fo$, $\Phi$, $\Phi^+$, $W$, $X_w$ and $C_w$ for $w\in W$ sinilarly to the finite-dimensional case, see, f.e., \cite{Kumar}. Our main result can be formulated as follows.

\mtheo{Let $G$ be the affine Kac--Moody group of type $\wt A_{n-1}$. Let $w_1$\textup,~$w_2$ be distinct involutions in the Weyl group $\wt S_n$ of $\Phi$. Then the tangent cones $C_{w_1}$ and~$C_{w_2}$ do~not coincide as subvarieties of the tangent space $T_p\Fo$.\label{mtheo:non_red}}

Two main technical tools we use are coadjoint orbits and embeddings of Weyl groups. Namely, one can naturally identify the tangent space $T_p\Fo$ with the dual space $\nt^*$ to the nilradical $\nt$ of the Lie algebra $\bt$ of the Borel subgroup $B$. To each involution $w\in W$ one can assign in a certain way the linear form $f_w\in\nt^*$ such that its $B$-orbit $\Omega_w$ under the coadjoint action of $B$ on $\nt^*$ is contained in~$C_w$, and $\dim\Omega_w=\dim C_w=l(w)$. On the other hand, consider an embedding of $G$ into a Kac--Moody group $\wt G$ such that $B\subseteq\wt B$ and $T\subseteq\wt T$, where $\wt T$ is a maximal torus of $\wt G$ contained in a Borel subgroup $\wt B$ of $\wt G$. Such an embedding induces the embeddings $w\subseteq\wt W$, $\Fo\subseteq\wt\Fo$ and $C_w\subseteq\wt C_w$, where $\wt\Fo=\wt G/\wt B$ and $\wt C_w$ is the tangent cone at the neutral point to the Schubert subvariety $\wt X_w$ of $\wt\Fo$ corresponding to an element $w\in W\subseteq\wt W$. Using coadjoint orbits mentioned above, we prove that if $w_1$ and $w_2$ are involutions in $W$ such that $C_{w_1}=C_{w_2}$ then $\wt l(w_1)=\wt l(w_2)$, where $\wt l$ is the length function on the Weyl group $\wt W$. Finally, given distinct involutions $w_1,w_2\in W$, we construct an embedding $G\subseteq\wt G$ such that $\wt l(w_1)\neq\wt l(w_2)$, which concludes the proof.

The paper is organized as follows. In Section~\ref{sect:Kac_Moody}, we briefly recall basic definitions and facts about Kac--Moody groups and algebras and about the structure of the group of affine permutations. In Section~\ref{sect:flags_cones} we discuss the geometry of flag ind-varieties for Kac--Moody groups and their tangent cones. We also describe the connections with the coadjoint action of the group $B$ and, based on this connection, present a proof of our main result modulo two rather technical results. These results are proved in Sections \ref{sect:proof1} and \ref{sect:proof2}.

\textsc{Acknowledgments}. We thank Evgeny Feigin for very fruitful discussions.

\sect{Basic definitions
}
\label{sect:Kac_Moody}
\addcontentsline{toc}{subsection}{\ref{sect:Kac_Moody}. Affine Kac--Moody group of type $\wt A_{n-1}$ and affine permutations}

In this section, we briefly introduce main definitions about an infinite-dimensional analogue of the symmetric group. Affine Weyl groups are obtained from usual Weyl groups by adding one additional generator. Weyl groups
and affine Weyl groups are special cases of Coxeter groups. Let $G$ be the set of invertible matrices
with entries in $\Cp(z)$, the field of rational functions of $z$ ($G$ is called a Kac--Moody group of type $\wt A_{n-1}$). The affine analog of the flag variety is $\Fo=G/B$, where $B$ is the set of all invertible matrices $g$ from $G$ with entries in $\Cp[[z]]$, the formal power series in $z$,
with the property that $\restr{g}{z=0}$ is upper triangular ($B$ is usually called the Iwahori subgroup of $G$). We use the same notation $X_w$ for Schubert varieties and affine Schubert
varieties, but the index comes from an affine Weyl group if $X_w$ is an affine Schubert variety.

Let \( S_n \) denote the permutation group, i.e., the group of all bijections taking \([n] = \{1, 2, \ldots, n\}\) to itself. We denote \( w \in S_n \) by its one-line notation \( w = [w_1, w_2, \ldots, w_n] \) if \( w \) maps \( i \) to \( w_i \). As a Coxeter group, \( S_n \) is generated by the adjacent transpositions \( s_i \) interchanging \( i \) with \( i + 1 \) with relations \( s_i^2 = e \), \( s_is_j = s_js_i \) for \(|i - j| \geq 2\), and \( s_is_{i+1}s_i = s_{i+1}s_is_{i+1} \) for all \( 1 \leq i \leq n - 2 \). Note that, clearly, \( S_n \) is the Weyl group of type \( A_{n-1} \) which corresponds to the complex reductive algebraic group \( G = \GL_n(\mathbb{C}) \).

Now, let us turn to the group of affine permutations. Let \( \wt{S}_n \) denote the set of all bijections \( w : \mathbb{Z} \to \mathbb{Z} \) with \( w(i + n) = w(i) + n \) for all \( i \in \mathbb{Z} \) and 
\begin{equation}\sum_{i=1}^n w(i) = \binom{n+1}{2}.\label{sum_w}
\end{equation}
Note that the permutation group $S_n$ can be naturally identified with the subgroup of $\wt S_n$, and, for $w\in S_n\subset\wt S_n$, condition (\ref{sum_w}) is satisfied automatically. The group \( \wt{S}_n \) is called the \emph{affine symmetric group}, and the elements of \( \wt{S}_n \) are called \emph{affine permutations}. This definition of affine permutations appeared in \cite[\S3.6]{Lusztig83}, see also \cite{BjornerBrenti}.

We can view an affine permutation in its one-line notation as the infinite string
$$\ldots w_{-1} w_0 w_1 w_2 \ldots w_n w_{n+1} \ldots,$$
where \( w_i = w(i) \). An affine permutation is completely determined by its action on any window of \( n \) consecutive indices \([w_i, w_{i+1}, \ldots, w_{i+n-1}]\). In particular, it is enough to record the base window \([w_1, \ldots, w_n]\) to capture all the information about \( w \). Sometimes, however, it will be useful to write down a larger section of the one-line notation or use a different window. For example, if $$w = [8, 1, 3, 5, 4, 0] \in \tilde{S}_6$$ in base window notation, then we could also say $$ w = [8, 1, 3, 5, 4, 0\mid14, 7, 9, 11, 10, 6, \ldots], $$ or \( w \) is determined by the window \([3, 5, 4, 0, 14, 7]\).

We can also view an affine permutation as a matrix. As mentioned before, the entries of the matrix will come from \( \mathbb{C}(z) \). For any \( w = [w_1, \ldots, w_n] \in \wt{S}_n \), write \( w_i = a_i + nb_i \), where \( 1 \leq a_i \leq n \) for each \( 1 \leq i \leq n \). Then \( \bar{w} := [a_1, \ldots, a_n] \in S_n \) and \( \sum b_i = 0 \). We will define the \textit{affine permutation matrix} for \( w \) to be the \( n \times n \) matrix \( M = (m_{ij}) \) with \( m_{\bar{w}_i, i} = z^{b_i} \) and all other entries equal 0. In the case where \( w = \bar{w} \in S_n \), this reduces to the usual notion of a permutation matrix.

Given \( i \neq j \) mod \( n \), let \( t_{ij} \) denote the affine transposition that interchanges \( i \) and \( j \) and leaves all \( k \) not congruent to \( i \) or \( j \) fixed. Since \( t_{ij} = t_{i+n,j+n} \), it suffices to assume \( 1 \leq i \leq n \) and \( i < j \). Note that if \( 1 \leq i < j \leq n \), the above notion of transposition is the same as for the symmetric group.

As a Coxeter group, \( \tilde{S}_n \) is the affine Weyl group of type \( A \). It is generated by
$$S = \{s_i := t_{i,i+1}: 1 \leq i \leq n - 1\} \cup \{s_0 := t_{n,n+1}\}.$$
We will denote the set of all affine transpositions by
$$ T = \{t_{ij}: 1 \leq i \leq n, \, i < j, \, i \not\equiv j \text{ mod } n\}.$$
The root system of $\wt{A}_{n-1}$ can be described in the following way: $\Phi=\Phi^+\cup\Phi^-$, $\Phi^-=-\Phi^+$ and $$\Phi^+=\{\alpha+k\delta\:|\: k\geq0\}\cup\{-\alpha+k\delta\:|\: k>0\}\cup\{k\delta\:|\: k>0\},$$ where $\alpha$ is a positive root in~$A_{n-1}$ and $\delta$ relates to $s_0$.

By definition, the \emph{length} $l(w)$ of an element $w\in\wt S_n$ is the length of a reduced expression of $w$ as a product of $s_i$'s (such an expression is called \emph{reduced} if it has minimal possible length). It is known that (see \cite[Lemma 4.2.2]{Shi86}), given $w\in\wt S_n$, one has
$$l(w) = \sum_{1 \le i < j \le n} \left|\left[ \dfrac{w_j - w_i}{n}\right] \right| = \operatorname{inv}(w_1, \ldots, w_n) + \sum_{1 \le i < j \le n} \left[ \dfrac{|w_j - w_i|}{n} \right],$$
where $\mathrm{inv}(w_1, \ldots, w_n) = |\{ 1 \le i < j \le n\colon w_i > w_j \}|$, and, as usual, $[x]$ denotes the maximal integer $z$ such that $z\le x$.

\exam{Let $w = [8, 1, 3, 5, 4, 0] \in \wt{S}_6$, as above. Then
\begin{equation*}
\begin{split}
&l(w)=\left( \left|\left[ \dfrac{1 - 8}{6}\right] \right|+\left|\left[ \dfrac{3- 8}{6}\right] \right|+\left|\left[ \dfrac{5-8}{6}\right] \right|+\left|\left[ \dfrac{4-8}{6}\right] \right|+\left|\left[ \dfrac{0 - 8}{6}\right] \right| \right) + \\
&\left( \left|\left[ \dfrac{3 - 1}{6}\right] \right|+\left|\left[ \dfrac{5 - 1}{6}\right] \right|+\left|\left[ \dfrac{4 - 1}{6}\right] \right|+\left|\left[ \dfrac{0 - 1}{6}\right] \right| \right)+   \left( \left|\left[ \dfrac{5 - 3}{6}\right] \right|+\left|\left[ \dfrac{4 - 3}{6}\right] \right|+\left|\left[ \dfrac{0 - 3}{6}\right] \right|\right) +\\    
&\left( \left|\left[ \dfrac{4-5}{6}\right] \right|+\left|\left[ \dfrac{0 - 5}{6}\right] \right|+ \right)+  \left|\left[ \dfrac{0-4}{6}\right] \right| =7+1+1+2+1=12.
\end{split}
\end{equation*}
One of its reduced expressions is the following:
$$w=s_1s_2s_3s_4s_5s_0s_4s_3s_2s_1s_4s_0.$$
}

Let $w\in \wt S_n$ be an involution, i.e., $w^2=\id$. The involution $w$ can be uniquely expressed as the product of pairwise commuting reflections: there exist unique distinct pairs of numbers $p_i,q_i$ belonging to the set $[1,n]=\{1,2,\ldots,n\}$, $i=1,\ldots,m$, such that $w(p_i)=q_i+k_in$ and $w(q_i)=p_i-k_in$ for some $k_i\in\Zp_{\geq0}$ (we assume that $p_i<q_i$ if $k_i=0$).

\defi{We call the set $\{(p_1, q_1, k_1),  (p_2, q_2, k_2), \dots,(p_m, q_m, k_m)\}$ the \emph{support} of $w$ and denote it by $\Supp{w}$.} 

\exam{i) The involution $s_0 \in \wt S_n$ is the reflection itself and can be expressed as \( s_0(1) = n - n\times1 \), \( s_0(2) = 2 \), \dots, \( s_0(n) = 1 +n\times1 \), hence $\Supp{s_0}=\{(1, n, 1)\}$.

ii) The involution $w=[6, -3, 4, 3]  \in \wt S_4$ can be uniquely expressed as the product of two pairwise commuting reflections, where the first one acts on $1, 2$: \( w(1) = 2 + 4\times1 \), \( w(2) = 1-4\times 1 \); while the second one acts on $3, 4$: \( w(4) = 3 + 4\times0 \), \( w(3) = 4-4\times 0 \). Thus, $\Supp{w}=\{(1, 2, 1), (4, 3, 0)\}$.
}

\sect{Sketch of the proof of the main result
}
\label{sect:flags_cones}
\addcontentsline{toc}{subsection}{\ref{sect:flags_cones}. Flag varieties, tangent cones and coadjoint orbits}

In this section we describe connections of the geometry of flag varieties for our Kac--Moody group and the orbits of the coadjoint action of the Iwahori group $B$ and give the proof of the Theorem~\ref{mtheo:non_red} modulo some technical results proved in the next sections. Denote by $\gt$, $\bt$, $\nt$ the Lie algebras of $G$, $B$, $U$ respectively. Namely, $\gt$ consists of all $n\times n$ matrices with entries in $\Cp(z)$, $\bt$ is the Lie algebra of matrices from $\gt$ with entries in $\Cp[[z]]$ which are upper-triangular if $z=0$, and $\nt$ is the Lie algebra of matrices from $\bt$ with zeroes on the diagonal for $z=0$. For each pair $(i,j)$ such that $1\leq i<j\leq n$ and each $k\in\Zp_{\geq0}$, one can consider the matrix $z^ke_{i,j}\in\nt$. We denote by $(z^ke_{i.j})^*$ the unique linear function on $\nt$ such that
\begin{equation*}
(z^ke_{i,j})^*(z^re_{p,q})=
\begin{cases}
1,&\text{if }i=p,j=q,k=r,\\
0&\text{otherwise}.
\end{cases}
\end{equation*}
We also denote by $\nt^-$ the Lie algebra of matrices $\lambda$ from $\gt$ such with entries in $\Cp[z^{-1}]$ such that all elements $\lambda_{i,j}$ for $i\leq j$ belong to $z^{-1}\Cp[z^{-1}]$. Clearly, $\gt=\bt\oplus\nt^-$. Furthermore, $\nt^-$ can be naturally identified with a subspace of $\nt^*$ via the trace from on $\gt$. Under this identification, $(z^ke_{i,j})^*$ corresponds to $z^{-k}e_{j,i}$.


The analogue of the flag variety, $\Fo=G/B$, can be naturally endowed with the structure of an ind-variety, see, e.g., \cite{Kumar02}. Hence, one can define the tangent space $T_p\Fo$, which is naturally isomorphic to the quotient vector space $\gt/\bt$. The group $B$ acts on $\Fo$ by conjugation. Since $p$ is $B$-stable, $B$ acts on the tangent space $T_p\Fo\cong\gt/\bt$, so $B$ acts on $\nt^-$. This action is called \emph{coadjoint}. We denote the result of coadjoint action by $b.\lambda$, $b\in B$, $\lambda\in\nt^-$. In 1962, for the finite-dimensional setting, A. Kirillov discovered that orbits of this action play an important role in representation theory of the group $B$ and of its unipotent radical, see, e.g., \cite{Kirillov62}, \cite{Kirillov04}.

Pick an involution $w\in\wt S_n$, let $\Supp{w}=\{(p_1,q_1, k_1),\ldots,(p_m,q_m, k_m)\}$ and $w(p_i)=q_i+k_in$, $w(q_i)=p_i-k_in$ for $k_i\in\Zp_{\geq0}$. Put $$f_w=\sum_{i=1}^m(z^{k_i}e_{p_i,q_i})^*\in\nt^-\subseteq\nt^*.$$

\defi{We say that the $B$-orbit $\Omega_w\subseteq\nt^-\subseteq\nt^*$ of $f_w$ is \emph{associated} with the involution $w$.}

One can easily check that $\Omega_w\subseteq C_w$. Indeed, suppose that $\Phi_0^+$ is a subset of $\Phi^+$ such that if $\alpha,\beta\in\Phi_0^+$ and $\alpha+\beta\in\Phi^+$, then $\alpha+\beta\in\Phi_0^+$. Then $\Phi_0=\Phi_0^+\cup(-\Phi_0^+)$ is a root subsystem of $\Phi$; denote by $G_0$ the corresponding Kac--Moody subgroup of $G$. Then $B_0=B\cap G_0$ is an Iwahori subgroup of $G_0$ contained in $B$, so there exists a natural embedding of the flag ind-varieties $\Fo_0=G_0/B_0\subseteq\Fo$. Furthermore, the Weyl group $W_0$ of $G_0$ can be realized as a subgroup of $W$ generated by $s_{\alpha}$, $\alpha\in\Phi_0$. Hence, if $w\in W_0$, then there exists a natural embedding of the Schubert varieties $X_{w,0}\subseteq X_w$ and of the tangent cones $C_{w,0}\subseteq C_w$. (Here $X_{w,0}$ is the Schubert subvariety of $\Fo_0$ corresponding to $w$, and $C_{w,0}$ is the tangent cone to~$X_{w,0}$ at the point~$p$.) Now, put $$\Phi_0^+=\bigcup_{(p_i,q_i, k_i)\in\Supp{w}}\{\epsi_{p_i}-\epsi_{q_i}+k_i\delta\}.$$ Since $\Phi_0$ is isomorphic as a root system to the direct product of $|\Supp{w}|$ copies of $A_1$ and $w$ is the longest element of $W_0$, $X_{w,0}$ coincides with $\Fo_0$, so $C_{w,0}$ is the entire tangent space $T_p\Fo_0$, i.e., the linear span of the root covectors $(z^{k_i}e_{p_i,q_i})^*$, $(p_i,q_i, k_i)\in\Supp{w}$. Hence $\Omega_w\subseteq C_{w,0}\subseteq C_w$, because $C_w$ is $B$-stable (in fact, the tangent cone to an arbitrary Schubert variety is $B$-stable).


Pick a linear form $\lambda\in\nt^-$ and, as usual, put
\begin{equation*}
\begin{split}
\Stab_B\lambda&=\{g\in B\mid g.\lambda=\lambda\}\subseteq B,\\
\Stab_{\bt}\lambda&=\{x\in\bt\mid\lambda([x,y])=0\text{ for all }y\in\nt\}\subseteq\bt.
\end{split}
\end{equation*}
We claim that 
\begin{equation}
\dim B.\lambda=\dim\bt/\Stab_{\bt}\lambda.\label{formula:dim_Omega_lambda}
\end{equation}
First, let $m\in\Zp_{\geq0}$ be the maximal degree of $\lambda_{i,j}$ as polynomials in $z^{-1}$. Then $(z^{m+1}B).\lambda=\lambda$ and $(z^{m+1}\bt).\lambda=0$, where
\begin{equation*}
\begin{split}
z^{m+1}\bt&=\{x\in\bt\mid\text{ all }x_{ij} \text{ are divisible by }z^{m+1}\},\\
z^{m+1}B&=\{g\in B\mid\text{ all }g_{ij} \text{ are divisible by }z^{m+1}\}.
\end{split}
\end{equation*}
Clearly, $z^{m+1}\bt$ is an ideal in $\bt$, while $z^{m+1}B$ is a normal subgroup in $B$.
Hence, only $B/\Stab_B\lambda$ and $\bt/(z^{m+1}\bt)$, which are finite-dimensional, could have non-trivial action on $\lambda$, thus, by standard algebraic-geometric arguments, we obtain the required equality.

Now, assume that $G''$ is a Kac--Moody group ot type $\wt A_{n-1}$, $G'$ is a Kac--Moody subgroup of $G''$, $T'$ (resp. $T''$) is a maximal torus of $G'$ (resp. of $G''$), $T'=T''\cap G'$, $B'$ (resp.~$B''$) is an Iwahori subgroup of $G'$ (resp. of $G''$) containing $T'$ (resp.~$T''$), $B'=B''\cap G'$, and $\Phi'$ (resp. $\Phi''$) is the root system of $G'$ (resp. of $G''$) with respect to $T'$ (resp. to~$T''$). We denote by $W'$ (resp. by $W''$) the Weyl group of $\Phi'$ (resp. of $\Phi''$). Denote by $\Fo'=G'/B'$, $\Fo''=G''/B''$ the flag ind-varieties. Put $p'=eB'\in\Fo'$, $p''=eB''\in\Fo''$. One can obviously define the Lie algebras $\gt'$,~$\bt'$,~$\nt'$ and $\gt''$, $\bt''$, $\nt''$ by the similar way as above. One can consider the space $\nt'^-\cong\gt'/\bt'$ as a~subspace of $\nt''^-\cong\gt''/\bt''$. Hence we can consider $T_{p'}\Fo'$ as a subspace of $T_{p''}\Fo''$.

Pick involutions $w_1,w_2\in W'$. Let $C_i'$ be the reduced tangent cone at the point $p'$ to the Schubert subvariety $X_{w_i}'$ of the flag variety $\Fo'$, $i=1,2$. Similarly, let $C_i''$ be the reduced tangent cone at $p''$ to the Schubert subvariety $X_{w_i}''$ of $\Fo''$, $i=1,2$. Denote by $l'$ (resp. by $l''$) the length function on the Weyl group $W'$ (resp. on $W''$). Assume $C_1'=C_2'$. This implies that $$l'(w_1)=l'(w_2).$$ Note that $C_i'\subseteq C_i''$, hence $B''.C_i'\subseteq C_i''$, $i=1,2$. Denote by $\Omega_{w_i}'\subseteq\nt'^*$ the coadjoint $B'$-orbit associated with the involution $w_i$, $i=1,2$; define $\Omega_{w_i}''$ by the similar way. It follows from Theorem \ref{theo:dim_Omega}  below that
\begin{equation*}
\begin{split}
l''(w_i)&=\dim C_i''\geq\dim B''.C_i'\geq\dim B''.\Omega_{w_i}'\\
&=\dim\Omega_{w_i}''=l''(w_i),
\end{split}
\end{equation*}
because $\Omega_{w_i}''=B''.\Omega_{w_i}'$. This implies $l''(w_i)=\dim C_i''=\dim B''.C_i'$. But $C_1'=C_2'$, thus $\dim C_1''=\dim C_2''$. We obtain the following result:
\begin{equation}
\text{if $C_1'=C_2'$, then $l''(w_1)=l''(w_2)$.}\label{formula:if_cones_then_ls}
\end{equation}

Finally, let $w_1$, $w_2$ be distinct involutions from the Weyl group $\wt S_n$ of affine permutations such that $C_{w_1}=C_{w_2}$ and, consequently, $l(w_1)=l(w_2)$. Then, by Proposition~\ref{lemm: embed}, there exists an embedding $\phi$ of $\wt S_n$ into $\wt S_{n+2}$ such that $l'(\phi(w_1))\neq l'(\phi(w_2))$\textup, where $l'$ denotes the length function on $\wt S_{n+2}$. This contradicts formula (\ref{formula:if_cones_then_ls}), which concludes the proof of our main result, Theorem~\ref{mtheo:non_red}.

\sect{The dimension of the orbit associated with an involution}
\label{sect:proof1}
\addcontentsline{toc}{subsection}{\ref{sect:proof1}. Proof}

In this section, we calculate the dimension of the orbit $\Omega_w$ associated with an involution $w$ from the group $\wt S_n$ of affine permutations. This is the first of the two key ingredients in the proof of our main result, see the previous section.

\theop{Let $w$ be an involution in the group of affine permutations $\wt S_n$. Then the dimension of the orbit $\Omega_w$ equals $l(w)$.\label{theo:dim_Omega}}{}Put $\lambda=f_w$. As we mentioned before (see formula (\ref{formula:dim_Omega_lambda}), $\dim \bt/\Stab_{\bt}\lambda = \dim\Omega_w$. It is clear that
\begin{equation*}
\Stab_{\bt}\lambda=\{x\in\bt\mid\lambda([x, z^k e_{ij}])=0\text{ for all } (i, j, k),\text{ where }k\geq 0\text{ and }i<j, \text{ or }k> 0\text{ and }i\geq j\}.
\end{equation*}
Hence, the codimension $\codim_{\bt}\Stab_{\bt}(\lambda)$ of the stabilizer in $\bt$ equals the number of linearly independent equations of the form $\lambda([x, z^k e_{ij}])=0$ for all such $i,j,k$. Clearly,
\begin{equation*}
\begin{split}
[x, z^k e_{ij}]&= z^k\left(\sum_{s=1}^nx_{si}e_{sj}-\sum_{s=1}^nx_{js}e_{is}\right) \\
&=z^k\begin{pmatrix}
0 & \cdots & 0 & x_{1i} & 0 & \cdots & 0 \\
\vdots& \ddots & \vdots & \vdots & \vdots & \ddots & \vdots \\
0 & \cdots & 0 & x_{i-1,i} & 0 & \cdots&0 \\
-x_{j1} & \cdots & -x_{j,j-1} & x_{ii} - x_{jj} & -x_{j,j+1} & \cdots & -x_{jn} \\
0&\cdots&0&x_{i+1,i}&0&\cdots&0\\
\vdots& \ddots & \vdots & \vdots & \vdots & \ddots & \vdots \\
0& \cdots & 0 & x_{ni} & 0 & \cdots & 0
\end{pmatrix}.
\end{split}
\end{equation*}
Note that this matrix can have nonzero elements only in the $i$th row and in the $j$th column.

Assume first that $\lambda=(z^se_{p, q})^*$, where $p\neq q$. We would like to compute the number of linearly independent equations mentioned above. There are four different types of equations depending on the indices $i$, $j$ of $e_{ij}$ (for $k=0$, $i$ should be less than $j$). $(x_{ij})_{s-k}$ defines the coefficient with $z^{s-k}$ of the polynomial (or row)  $x_{ij}(z)$.
\begin{enumerate}
    \item $i=p$, $j=q$: $\lambda([x, z^k e_{ij}])=(x_{pp}-x_{qq})_{s-k}=0$, one equation for each $k$, where
    \begin{equation*}
    k\in\begin{cases}
    [0,s],&\text{if }p<q,\\
    [1,s],&\text{if }p>q.\\
    \end{cases}
    \end{equation*}
    
    \item $i =j= p$, or $i=j=q$ (they give the same equations): $\lambda([x, z^k e_{ij}])=(x_{pq})_{s-k}=0$,  where \begin{equation*}
    k\in\begin{cases}
    [1,s],&\text{if } p<q,\\
    [1,s-1],&\text{if } p>q.\\
    \end{cases}
    \end{equation*}
    
    \item  $i=p$, $j\neq q, p$: $\lambda([x, z^k e_{ij}])=(x_{jq})_{s-k}=0$, where
    \begin{equation*}
    k\in\begin{cases}
    [0,s],&\text{if }i<j<q,\\
    [1,s],&\text{if }i> j,~j<q,\\
    [0,s-1],&\text{if }i<j, \: j>q,\\
    [1,s-1],&\text{if }i> j>q.\\
    \end{cases}
    \end{equation*}
    
    \item  $i \neq p, q$, $j=q$: $\lambda([x, z^k e_{ij}])=(x_{pi})_{s-k}=0$,  where \begin{equation*}
    k\in\begin{cases}
    [0,s],&\text{if }p<i<j,\\
    [1,s],&\text{if }i>j,~p<i,\\
    [0,s-1],&\text{if }i<j, \: p>i,\\
    [1,s-1],&\text{if }i> j, \: p>i.\\
    \end{cases}
    \end{equation*}
    \item otherwise, $\lambda([x, z^k e_{ij}])=0_{s-k}=0$ (trivial equations).
\end{enumerate}



It is well known that $l(w)$ equals the number of positive roots $\beta\in \Phi^+$ such that $w(\beta)\in \Phi^-$. For $\wt{A}_{n-1}$,  $\Phi^+=\{\alpha+k\delta\:|\: k\geq0\}\cup\{-\alpha+k\delta\:|\: k>0\}\cup\{k\delta\:|\: k>0\}$, where $\alpha$ is a positive root in~$A_{n-1}$. Because imaginary positive roots of the form $k\delta$, $k>0$, can not become negative after applying $w$, we can consider only real positive roots, i.e., the roots from the set $$\Phi^+_r=\{\alpha+k\delta\:|\: k\geq0\}\cup\{-\alpha+k\delta\:|\: k>0\}.$$
The support of the reflection $s_{\alpha}$ with respect to a positive root $\alpha=e_p-e_q+s\delta$ has the form $\Supp{s_{\alpha}}=\{(p, q, s)\}$. We will denote $s_{\alpha}$ by $w_{(p, q, s)}$.

Pick a positive root $\beta=e_i-e_j+k\delta$, with $i\neq j$, $k\geq0$ and $k=0$ only if $i<j$, then
\begin{equation*}
\begin{split}
&s_\alpha(\beta)=e_i-e_j+k\delta  - (\delta_{ip}+\delta_{jq}-\delta_{iq}-\delta_{jp})(e_p-e_q+s\delta).
\end{split}
\end{equation*}
Now we would like to know when $s_\alpha(\beta)$ is negative for different cases of triples $(i, j, k)$.
\begin{enumerate}
    \item  $i=p$, $j=q$: $s_\alpha(\beta)=e_p-e_q+k\delta  - 2(e_p-e_q+s\delta)=e_q-e_p+(k-2s)\delta$; it will be negative when \begin{equation*}
    k\in\begin{cases}
    [1,2s-1],&\text{if }p>q,\\
    [0,2s],&\text{if }p<q,\\
    \end{cases}
    \end{equation*}
    \item  $i=q$, $j=p$: $s_\alpha(\beta)=e_q-e_p+k\delta  + 2(e_p-e_q+s\delta)=e_p-e_q+(k+2s)\delta$, it  could be negative if $k=s=0$. Then $q<p$ because of positivity of $\beta$ and, at the same time, $p<q$ because of positivity of $\alpha$, a contradiction.
    \item  $i=p$, $j\neq q$: $s_\alpha(\beta)=e_p-e_j+k\delta  -(e_p-e_q+s\delta)=e_q-e_j+(k-s)\delta$; it will be negative when \begin{equation*}
    k\in\begin{cases}
    [0,s],&\text{if }i<j<q,\\
    [1,s],&\text{if }i> j,~j<q,\\
    [0,s-1],&\text{if }i<j, \: j>q,\\
    [1,s-1],&\text{if }i> j>q.\\
    \end{cases}
    \end{equation*} 
    \item $i\neq p$, $j=q$: $s_\alpha(\beta)=e_i-e_q+k\delta  -(e_p-e_q+s\delta)=e_i-e_p+(k-s)\delta$; if will be negative when\begin{equation*}
    k\in\begin{cases}
    [0,s],&\text{if }p<i<j,\\
    [1,s],&\text{if }i>j,~p<i,\\
    [0,s-1],&\text{if }i<j, \: p>i,\\
    [1,s-1],&\text{if }i> j, \: p>i.\\
    \end{cases}
    \end{equation*}
    \item  $\{i, j\}\cap\{p,q\}=\varnothing$: $s_\alpha(\beta)=\beta$ is positive.
    \item $i=q$, $j\neq p$: $s_\alpha(\beta)=e_q-e_j+k\delta  +(e_p-e_q+s\delta)=e_p-e_j+(k+s)\delta$, it  could be negative if $k=s=0$ then $q<j$ because of positivity of $\beta$, and $p<q$ because of positivity of $\alpha$, but then $s_\alpha(\beta)$ is positive.
    \item $i\neq q$, $j=p$: $s_\alpha(\beta)=e_i-e_p+k\delta  +(e_p-e_q+s\delta)=e_i-e_q+(k+s)\delta$, it  could be negative if $k=s=0$ then $i<p$ because of positivity of $\beta$, and $p<q$ because of positivity of $\alpha$, but then $s_\alpha(\beta)$ is positive. 
\end{enumerate}
Then $\beta$ becomes negative in cases $1$, $3$, $4$. 

Let us construct a bijection between the set $\Phi^\pm_{w_{(p, q, s)}}$ of positive roots becoming negative after applying $w_{(p, q, s)}$, and linearly independent equations given by $\lambda=(z^se_{pq})^*$. To each $\beta \in \Phi^\pm_{w_{(p, q, s)}}$  associated with triple $(i,j,k)$ and related to cases 3 and 4, we attach the equations of third and fourth types respectively:  $\lambda([X, z^k e_{ij}])=0$. If $\beta$ is from case 1 (i.e., $i=p$, $j=q$), then there are two opportunities. \begin{enumerate}

\item $p>q$:
\begin{equation*}
    k\in\begin{cases}
    [1,s],&\text{we attach to $\beta$ the equation  }\lambda([x, z^k e_{pq}])=(x_{pp}-x_{qq})_{s-k}=0,\\
    [s+1,2s-1],&\text{we attach to $\beta$ the equation  }\lambda([x, z^k e_{pp}])=(x_{pq})_{2s-k}=0.\\
    \end{cases}
    \end{equation*} . 
\item $p<q$: \begin{equation*}
k\in\begin{cases}
    [0,s],&\text{we attach to $\beta$ the equation }\lambda([x, z^k e_{pq}])=(x_{pp}-x_{qq})_{s-k}=0,\\
    [s+1,2s],&\text{we attach to $\beta$ the equation  }\lambda([x, z^k e_{pp}])=(x_{pq})_{2s-k}=0.\\
    \end{cases}
    \end{equation*}
\end{enumerate}

\noindent Thus, we have proved the theorem in the case when $w$ is just a reflection $s_{\alpha}$.

Now, consider an arbitrary involution $w= w_{(p_1,q_1, k_1)}w_{(p_2,q_2, k_2)}\ldots w_{(p_m,q_m, k_m)}$. We claim that, given an entry $(i,j)$, there are at most two triples in the support of $w$ that can affect the equations involving $x_{i,j}$. Indeed, $\lambda=\sum_{i=1}^m\lambda_{w_{(p_i,q_i,k_i)}}$ and $x_{i,j}$ can occur in an equation corresponding to $\lambda_{w_{(p,q,k)}}$ only if $i=p$ or $j=q$. Similarly, given a positive root $\beta$, there are at most two triples in the support of $w$ that can affect the negativity or positivity of $w(\beta).$

Consider all positive roots $\beta=e_i-e_j+k\delta$ affected by triples $(a,b,l)$ and $(p,q,s)$, such that $i \in \{p, q\}$ and $j \in \{a, b\}$. Compute the number of them subsequently. 
\begin{enumerate}
\item $\beta=e_p-e_b+k\delta$:
    \begin{equation*}
        w'(\beta)=e_q-e_a+(k-l-s)\delta,\\
    \end{equation*}    
    \begin{equation*}
        k\in \begin{cases}
        [0,l+s],&\text{if } q> a, \: p < b, \\
        [1,l+s],&\text{if } q> a, \: p > b,\\
        [0,l+s-1],&\text{if } q<a, \: p < b, \\
        [1,l+s-1],&\text{if } q< a, \: p > b. \\
        \end{cases}
    \end{equation*}

\item $\beta=e_p-e_a+k\delta$:
\begin{equation*}
    w'(\beta)=e_q-e_b+(k+l-s)\delta,
\end{equation*}
 \begin{equation*}
        k\in \begin{cases}
        [0,s-l],&\text{if } s\geq l, \: p < a, \: q>b, \\
        [1,s-l],&\text{if }s\geq l, \: p > a, \: q>b,\\
        [0,s-l-1],&\text{if } s\geq l, \: p < a, \: q<b, \\
        [1,s-l-1],&\text{if } s\geq l, \: p > a, \: q<b, \\
        \varnothing ,&\text{if } s< l. \\
        \end{cases}
    \end{equation*}

\item $\beta=e_q-e_a+k\delta$:
\begin{equation*}
    w'(\beta)=e_p-e_b+(k+l+s)\delta,
\end{equation*}
If $w'(\beta)$ is negative, then  $k=l=s=0$ and $p>b$, therefore $q<a$, $a<b$, $p<q$, but it is impossible.

\item $\beta=e_q-e_b+k\delta$:
\begin{equation*}
    w'(\beta)=e_p-e_a+(k+s-l)\delta,
\end{equation*}
 \begin{equation*}
        k\in \begin{cases}
        [0,s-l],&\text{if } l\geq s, \: p > a, \: q<b, \\
        [1,s-l],&\text{if }l\geq s, \: p > a, \: q>b,\\
        [0,l-s-1],&\text{if } l\geq s, \: p < a, \: q<b, \\
        [1,l-s-1],&\text{if } l\geq s, \: p < a, \: q>b, \\
        \varnothing ,&\text{if } l< s. \\
        \end{cases}
    \end{equation*}
\end{enumerate}
Consider all equations $\lambda([x, z^k e_{ij}])=0$, such that $i \in \{p, q\}$ and $j \in \{a, b\}$. Compute the number of them subsequently. 
\begin{enumerate}
\item $\lambda([x, z^k e_{pb}])=(x_{bq})_{s-k}+(x_{ap})_{l-k}=0$:
    \begin{equation*}
          k\in \begin{cases}
        [0,l],&\text{if } l> s, \: p < b, \: a < p, \\
        [1,l],&\text{if } l> s, \: p > b, \:  a<p,\\
        [0,l-1],&\text{if } l> s, \: p < b, \: a > p, \\
        [1,l-1],&\text{if } l> s, \: p > b, \:  a>p,\\
        [0,s],&\text{if } l< s, \: p < b, \: b<q, \\
        [1,s],&\text{if } l< s, \: p > b, \: b<q, \\
        [0,s-1],&\text{if } l< s, \: p < b, \: b>q, \\
        [1,s-1],&\text{if } l< s, \: p > b, \: b>q, \\
        [0,s],&\text{if } l= s, \: p < b, \: (b<q \textup{ or } a < p),\\
        [1,s],&\text{if } l=s, \: p > b, \: (b<q\textup{ or } a < p), \\
        [0,s-1],&\text{if } l= s, \: p < b, \: (b>q\textup{ and } a > p), \\
        [1,s-1],&\text{if } l= s, \: p > b, \: (b>q\textup{ and } a > p). \\
        \end{cases}
    \end{equation*}

\item $\lambda([x, z^k e_{pa}])=(x_{aq})_{s-k}=0$ and $\lambda([x, z^k e_{qb}])=(x_{aq})_{l-k}=0$:

 \begin{equation*}
        k\in \begin{cases}
        [0,s],&\text{if } s> l, \: p < a , \: a<q, \\
        [1,s]*,&\text{if }s> l, \: p > a, \: a<q,\\
        [0,s-1],&\text{if } s> l, \: p < a, \: q<a, \\
        [1,s-1]*,&\text{if } s> l, \: p > a, \: q<a, \\
        [0,l],&\text{if } l> s, \: q > a,\: q<b, \\
        [1,l]*,&\text{if }l> s, \: q > a, \: q>b,\\
        [0,l-1],&\text{if } l> s, \: q < a, \: q<b, \\
        [1,l-1]*,&\text{if } l>s, \: q < a, \: q>b, \\
        [0,s],&\text{if } s= l, \: (p < a \textup{ or } q<b), \: a<q, \\
        [1,s],&\text{if }s=l, \: (p > a \textup{ and } q>b), \: a<q,\\
        [0,s-1],&\text{if } s= l, \: (p < a \textup{ or } q<b), \: q<a, \\
        [1,s-1],&\text{if } s= l, \: (p > a \textup{ and } q<b), \: q<a. \\  
        
        \end{cases}
    \end{equation*}
Here $*$ means that another equation can give us an option $k=0$, but would not give new linearly independent equation.

\item  $\lambda([x, z^k e_{qa}])=0=0$:
The equation is trivial. 
\end{enumerate}

Here, case 1 (respectively, 2, 3) corresponds to case 1 (respectively, 2 and 4, 3) for positive roots becoming negative. All these equations are linearly independent (except for two types of equations in case 2), and they are also linearly independent with other equations (because all other equations have different indices of $x_{ij}$).

Assume without loss of generality that $l\geq s$ (we can do so because of the symmetry). Firstly, assume $l> s$ and calculate the amount $N$ of the equations in cases 1 and 2 above. If $p<b$, $a<p$, $q>a$, $q<b$, then $N=2l+2$ linearly independent equations occur. For simplicity, given an inequality~$I$, we will write
\begin{equation*}
(I)=\begin{cases}1,&\text{if $I$ is true},\\
0,&\text{if $I$ is false}.
\end{cases}
\end{equation*}
Note that each inequality $p<b$, $a<p$, $q>a$, $q<b$ give us an extra equation (pushing left or right boundaries of $k$ by one), thus, $$N=2l-2+(p<b)+ (a<p) +(q>a)+ (q<b),$$ and this is the same amount that the number of positive roots becoming negative calculated previously. 

Secondly, if $l=s$, then $$N=2s-2 + (q>a) + (p<b) + (p<a \textup{ or } q<b)+ (p>a \textup{ or } q>b),$$ while the number of positive roots becoming negative equals $$2s-1 + (q>a) + (p<b) + (p<a \textup{ and } q>b)+ (p>a \textup{ and }
q<b).$$ Below we show that these values are the same:
\begin{equation*}
\begin{split}
N&=2s-2 + (q>a) + (p<b) + (p<a \textup{ or } q<b)+ (p>a \textup{ or } q>b)\\
&=2s-1 + (q>a) + (p<b) +(p<a \textup{ or } q<b)+ (p>a \textup{ or } q>b)-1\\
&=2s-1 + (q>a) + (p<b) + (p<a \textup{ and } q>b)+ (p>a \textup{ and }
q<b).
\end{split}
\end{equation*}

We see that the number of positive roots becoming negative because of two triples equals the number of linearly independent equations on $\Stab_{\bt}\lambda$ associated with these two triples. As we checked before, the same is true for the case of one triple. Thus, the result follows.
\newpage
\sect{Embeddings of groups of affine permutations}
\label{sect:proof2}
\addcontentsline{toc}{subsection}{\ref{sect:proof2}. Proof}

In this section, we construct an embedding of $\wt S_n$ to $\wt S_{n+2}$, which sends distinct involutions in $\wt S_n$ to permutations with different lengths. This is the second of the two key ingredients in the proof of our main result, see Section~\ref{sect:flags_cones}.

\propp{Let \label{lemm: embed} $w_1$, $w_2$ be distinct involutions from the Weyl group $\wt S_n$ of affine permutations such that $C_{w_1}=C_{w_2}$ \textup(and\textup, consequently\textup, $l(w_1)=l(w_2)$\textup). Then there exists an embedding $\phi$ of $\wt S_n$ into $\wt S_{n+2}$ such that $l'(\phi(w_1))\neq l'(\phi(w_2))$\textup, where $l'$ denotes the length function on $\wt S_{n+2}$.}
{Consider the embedding $\phi$ of $\wt S_n$ into $\wt S_{n+2}$, "adding" a pair of elements $0\leq a < c \leq n +1$ and a parameter $p$ that is for every $w\in\wt S_n$, $\phi(w)=w'$ and:
\begin{enumerate}
    \item if $-1<x<a$: $w'(x)=w(x)+ \left[\dfrac{w(x)}{n}\right](n+2)$,
    \item if $x=a$: $w'(x)=c+p(n+2)$,
    \item if $a<x<c$: $w'(x)=w(x-1)+ \left[\dfrac{w(x-1)}{n}\right](n+2)$,
    \item if $x=c$: $w'(x)=a-p(n+2)$,
    \item if $c<x<n+2$: $w'(x)=w(x-2)+ \left[\dfrac{w(x-2)}{n}\right](n+2)$.
\end{enumerate}

Consider $A=\{ 0, 1, \dots, a-1\}$,  $B=\{ a+1, a+2, \dots, c-1\}$ and  $C=\{ c+1, c+2, \dots, n+1\}$. 
Then
\begin{equation*}
\begin{split}
l'(w')&=\sum_{0\leq i<j<n+2}\left| \left[ \frac{w'(j)-w'(i)}{n+2}\right]\right|\\
&=l(w)+\left| \left[ \frac{w'(c)-w'(a)}{n+2}\right]\right|+\sum\limits_{i \in A}\left(\left| \left[ \frac{w'(a)-w'(i)}{n+2}\right]\right|+\left| \left[ \frac{w'(c)-w'(i)}{n+2}\right]\right|\right)\\
&+ \sum\limits_{i \in B}\left(\left| \left[ \frac{w'(i)-w'(a)}{n+2}\right]\right|+\left| \left[ \frac{w'(c)-w'(i)}{n+2}\right]\right|\right)\\
&+ \sum\limits_{i \in C}\left(\left| \left[ \frac{w'(i)-w'(a)}{n+2}\right]\right|+\left| \left[ \frac{w'(i)-w'(c)}{n+2}\right]\right|\right)
\end{split}
\end{equation*}
Assume $w'(i)=j+k(n+2)$, $j \in [0, n+1]$ and $p\geq 0$. Then consider several cases to understand what contribution the element $i$ makes to the last three terms as it goes to $j+k(n+2)$.
\begin{enumerate}
    \item if $i \in A, j\in A$; 
    \begin{equation*}
    \begin{split}
    \left| \left[ \frac{w'(a)-w'(i)}{n+2}\right]\right|+\left| \left[ \frac{w'(c)-w'(i)}{n+2}\right]\right|=|p-k| +|p+k|=2\max(|p|, |k|),
    \end{split}
    \end{equation*}
    
    \item if $i \in C, j\in C$,
    \begin{equation*}
    \begin{split}
    \left| \left[ \frac{w'(i)-w'(a)}{n+2}\right]\right|+\left| \left[ \frac{w'(i)-w'(c)}{n+2}\right]\right|=|p-k| +|p+k|=2\max(|p|, |k|),
    \end{split}
    \end{equation*}
    
    \item if $i \in A, j\in B$, if $k+p\geq0$:
    \begin{equation*}
    \begin{split}
    &\left| \left[ \frac{w'(a)-w'(i)}{n+2}\right]\right|+\left| \left[ \frac{w'(c)-w'(i)}{n+2}\right]\right|\\
    &=\left| \left[ \frac{c+p(n+2)-j-k(n+2)}{n+2}\right]\right|+\left| \left[ \frac{a-p(n+2)-j-k(n+2)}{n+2}\right]\right|=|p-k| +|p+k+1|\\
    &=2\max(|p|, |k|)+1,
    \end{split}
    \end{equation*}
    else if $k+p<0$:
    \begin{equation*}
    \begin{split}
    &\left| \left[ \frac{w'(a)-w'(i)}{n+2}\right]\right|+\left| \left[ \frac{w'(c)-w'(i)}{n+2}\right]\right|\\
    &=\left| \left[ \frac{c+p(n+2)-j-k(n+2)}{n+2}\right]\right|+\left| \left[ \frac{a-p(n+2)-j-k(n+2)}{n+2}\right]\right|=|p-k| +|p+k+1|\\
    &=2\max(|p|, |k|)-1,
    \end{split}
    \end{equation*}
    
    \item if $i \in B, j\in A$ (because of symmetry it is identical to previous case with $-k$), if $p- k\geq0$:
      \begin{equation*}
    \begin{split}
    &\left| \left[ \frac{w'(i)-w'(a)}{n+2}\right]\right|+\left| \left[ \frac{w'(c)-w'(i)}{n+2}\right]\right|=2\max(|p|, |k|)+1,
    \end{split}
    \end{equation*}
    else if $p-k<0$:
    \begin{equation*}
    \begin{split}
    &\left| \left[ \frac{w'(a)-w'(i)}{n+2}\right]\right|+\left| \left[ \frac{w'(c)-w'(i)}{n+2}\right]\right|=2\max(|p|, |k|)-1,
    \end{split}
    \end{equation*}
    
    \item if $i \in A, j\in C$,
    \begin{equation*}
    \begin{split}
    &\left| \left[ \frac{w'(a)-w'(i)}{n+2}\right]\right|+\left| \left[ \frac{w'(c)-w'(i)}{n+2}\right]\right|\\
    &=\left| \left[ \frac{c+p(n+2)-j-k(n+2)}{n+2}\right]\right|+\left| \left[ \frac{a-p(n+2)-j-k(n+2)}{n+2}\right]\right|=|p-k-1| +|p+k+1|\\
    &=2\max(|p|, |k+1|),
    \end{split}
    \end{equation*}
    
    \item if $i \in C, j\in A$ (because of symmetry it is identical to previous case with $-k$),
    \begin{equation*}
    \begin{split}
    &\left| \left[ \frac{w'(i)-w'(a)}{n+2}\right]\right|+\left| \left[ \frac{w'(i)-w'(a)}{n+2}\right]\right|=2\max(|p|, |k-1|),
    \end{split}
    \end{equation*}
    
    \item if $i \in B, j\in B$,
    \begin{equation*}
    \begin{split}
    &\left| \left[ \frac{w'(i)-w'(a)}{n+2}\right]\right|+\left| \left[ \frac{w'(c)-w'(i)}{n+2}\right]\right|\\
    &=\left| \left[ \frac{j+k(n+2)-c-p(n+2)}{n+2}\right]\right|+\left| \left[ \frac{a-p(n+2)-j-k(n+2)}{n+2}\right]\right|=|k-p-1| +|k+p+1|\\
    &=2\max(|p+1|, |k|),
    \end{split}
    \end{equation*}
    
    \item if $i \in B, j\in C$,if $k\geq -p$:
    \begin{equation*}
    \begin{split}
     &\left| \left[ \frac{w'(i)-w'(a)}{n+2}\right]\right|+\left| \left[ \frac{w'(c)-w'(i)}{n+2}\right]\right|\\
    &=\left| \left[ \frac{j+k(n+2)-c-p(n+2)}{n+2}\right]\right|+\left| \left[ \frac{a-p(n+2)-j-k(n+2)}{n+2}\right]\right|=|k-p| +|k+p+1|\\
    &=2\max(|p|, |k|)+1,
    \end{split}
    \end{equation*}
    else if $k<-p$:
    \begin{equation*}
    \begin{split}
     &\left| \left[ \frac{w'(i)-w'(a)}{n+2}\right]\right|+\left| \left[ \frac{w'(c)-w'(i)}{n+2}\right]\right|\\
    &=\left| \left[ \frac{j+k(n+2)-c-p(n+2)}{n+2}\right]\right|+\left| \left[ \frac{a-p(n+2)-j-k(n+2)}{n+2}\right]\right|=|k-p| +|k+p+1|\\
    &=2\max(|p|, |k|)-1,
    \end{split}
    \end{equation*}

    \item if $i \in C, j\in B$ (because of symmetry it is identical to previous case with $-k$), if $p- k\geq0$:
    \begin{equation*}
    \begin{split}
     &\left| \left[ \frac{w'(i)-w'(a)}{n+2}\right]\right|+\left| \left[ \frac{w'(i)-w'(c)}{n+2}\right]\right|=2\max(|p|, |k|)+1,
    \end{split}
    \end{equation*}
    else if $p-k<0$:
    \begin{equation*}
    \begin{split}
     &\left| \left[ \frac{w'(i)-w'(a)}{n+2}\right]\right|+\left| \left[ \frac{w'(i)-w'(c)}{n+2}\right]\right|=2\max(|p|, |k|)-1.
    \end{split}
    \end{equation*}

\end{enumerate}
Set the $w'$-order of the number $x$ from $S=A\cup B\cup C$ to be the number $k$ such as $w'(x)=y+k(n+2)$ and $y \in [0, n+1]$.  Since we are only considering involutions, for each number $x$ with $w'$-order $k\neq0$, the number $y$ has $w'$-order $-k$. 
Denote the set of all the numbers with $w'$-order equal to $k$ by $V^k$. Consider $p=m$, where $m$ is the maximum number such that $V^m$ is non-empty. Then calculate values for the previous cases:
\begin{enumerate}
    \item if $i \in A, j\in A$: $2p$,
    \item if $i \in C, j\in C$: $2p$,
    \item if $i \in A, j\in B$: $2p+1$,
    \item if $i \in B, j\in A$: $2p+1$,
    \item if $i \in A, j\in C$: for $i \in V^m$, $2p+2$, else $2p$,
    \item if $i \in C, j\in A$: for $i \in V^{-m}$, $2p+2$, else $2p$,
    \item if $i \in B, j\in B$: $2p+2$,
    \item if $i \in B, j\in C$: $2p+1$,
    \item if $i \in C, j\in B$: $2p+1$.
\end{enumerate}

Now we can simplify our length formula:
\begin{equation*}
\begin{split}
l'(w')&=l(w)+\left| \left[ \frac{w'(c)-w'(a)}{n+2}\right]\right|+2pn+2|B|\\
&+4\#\{x \in A\cap V^m\mid w'(x)=y+m(n+2), y\in C\}.
\end{split}
\end{equation*}
Then
\begin{equation*}
\begin{split}
l'(w'_1)-l'(w'_2)
&=4(\#\{x \in A_1\cap V_1^m\mid w'_1(x)=y+m(n+2), y\in C_1\}\\
&-\#\{x \in A_2\cap V^m_2\mid w'_2(x)=y+m(n+2), y\in C_2\}).
\end{split}
\end{equation*}
\\
Assume that $V_1^m\neq V_2^m$, or  $V_1^m= V_2^m$  and $\restr{w_1}{V_1^m}\neq \restr{w_2}{V_2^m}$. Then choose the smallest numbers $t_1\in V_1^m$, $t_2\in V_2^m$ such that $w_1'(t_1)\neq w_2'(t_1)$, $w_1'(t_2)\neq w_2'(t_2)$, and assume $t_1\leq t_2$. If $t_1\neq t_2$, then put $a=t_1+1$, $c=t_1+2$ (we did not use the fixed $a, c$ in the proof before, so we can fix them now). Then \begin{equation*}
\begin{split}
l'(w'_1)-l'(w'_2)
&=4(\#\{x \in A_1\cap V_1^m\mid w'_1(x)=y+m(n+2), y\in C_1\}\\
&-\#\{x \in A_2\cap V^m_2\mid w'_2(x)=y+m(n+2), y\in C_2\})=4.
\end{split}
\end{equation*}
If $t_1=t_2=t$, then assume $l_1=w'_1(t)<w_2(t)=l_2$, and consider $a=t+1$, $c=l_1+1$:
\begin{equation*}
\begin{split}
l'(w'_1)-l'(w'_2) &=4(\#\{x \in A_1\cap V_1^m\mid w'_1(x)=y+m(n+2), y\in C_1\}\\
&-\#\{x \in A_2\cap V^m_2\mid w'_2(x)=y+m(n+2), y\in C_2\})=4
\end{split}
\end{equation*}
This is why $V_1^m= V_2^m$ and $\restr{w_1}{V_1^m}=\restr{w_2}{V_2^m}$, and therefore $V_1^{-m}= V_2^{-m}$ with $\restr{w_1}{V_1^{-m}}=\restr{w_2}{V_2^{-m}}$.
We will use induction to continue the proof. Assume $V_1^q\neq V_2^q$, $q\geq0$ and for each $q'>q$ one has $V_1^{q'}= V_2^{q'}$ and $\restr{w_1}{V_1^{q'}}=\restr{w_2}{V_2^{q'}}$. Define $P=\bigcup_{q'>q}(V^{q'}\cup V^{-q'})$ and $R = S\backslash P$. Then 
\begin{equation*}
\begin{split}
l'(w')&=l(w)+\left| \left[ \frac{w'(c)-w'(a)}{n+2}\right]\right| +\sum\limits_{i \in P}\left(\left| \left[ \frac{(w'(i)-w'(a))\frac{i-a}{|i-a|}}{n+2}\right]\right|+\left| \left[ \frac{(w'(c)-w'(i))\frac{c-i}{|c-i|}}{n+2}\right]\right|\right)\\
&=l(w)+\left| \left[ \frac{w'(c)-w'(a)}{n+2}\right]\right|+2p|R|+2|B\cap R|+4\#\{x \in A\cap V^q \mid w'(x)=y+q(n+2), y\in C\},
\end{split}
\end{equation*}
\begin{equation*}
\begin{split}
l'(w'_1)-l'(w'_2) &=4(\#\{x \in A_1\cap V_1^q\mid w'_1(x)=y+q(n+2), y\in C_1\}\\
&-\#\{x \in A_2\cap V^q_2\mid w'_2(x)=y+q(n+2), y\in C_2\})
\end{split}
\end{equation*}
because of the inductive assumption.
 
Now we will use the line of arguing similar to the case of $V^m$. Namely, assume $V_1^q\neq V_2^q$, or $V_1^q= V_2^q$ and $\restr{w_1}{V_1^q}\neq\restr{w_2}{V_2^q}$). Then choose the smallest numbers $t_1\in V_1^q$, $t_2\in V_2^q$  such that $w_1'(t_1)\neq w_2'(t_1)$, $w_1'(t_2)\neq w_2'(t_2)$, and assume $t_1\leq t_2$.

If $t_1\neq t_2$, consider that $a=t_1+1$, $c=t_1+2$, $p=q$ (we did not fix $a$, $c$ and $p$ in the proof before, so we can fix them now). Then \begin{equation*}
\begin{split}
l'(w'_1)-l'(w'_2)&=4(\#\{x \in A_1\cap V_1^q\mid w'_1(x)=y+q(n+2), y\in C_1\}\\
&-\#\{x \in A_2\cap V^q_2\mid w'_2(x)=y+q(n+2), y\in C_2\})
=4.
\end{split}
\end{equation*}
Else $t_1=t_2=t$, assume $l_1=w'_1(t)<w_2(t)=l_2$, then consider $a=t+1$, $c=l_1+1$:
\begin{equation*}
\begin{split}
l'(w'_1)-l'(w'_2)&=4(\#\{x \in A_1\cap V_1^q\mid w'_1(x)=y+m(n+2), y\in C_1\}\\
&-\#\{x \in A_2\cap V^q_2\mid w'_2(x)=y+q(n+2), y\in C_2\})=4.
\end{split}
\end{equation*}
This is why $V_1^q= V_2^q$ and $\restr{w_1}{V_1^q}=\restr{w_2}{V_2^q}$ for all $q$. Thus, $w_1=w_2$, a contradiction. This concludes the proof.}


\medskip\textsc{Stepan Bondar: National Research University Higher School of Economics,\break\indent ul. Usacheva 6, 119048, Moscow, Russia}

\emph{E-mail address}: \texttt{mstbondar@gmail.com}

\medskip\textsc{Mikhail Ignatev: National Research University Higher School of Economics,\break\indent Pokrovsky Boulevard 11, 109028, Moscow, Russia}

\emph{E-mail address}: \texttt{mihail.ignatev@gmail.com}


\begin{thebibliography}{XXXXXX}\addcontentsline{toc}{subsection}{References}


\bibitem[BB05]{BjornerBrenti} A. Bjorner, F. Brenti. Combinatorics of Coxeter groups. Graduate Texts in Mathematics \textbf{231}, Springer, 2005.

\bibitem[BIS16]{BochkarevIgnatyevShevchenko} M.A. Bochkarev, M.V. Ignatyev, A.A. Shevchenko. Tangent cones to Schubert varieties in types $A_n$, $B_n$ and $C_n$. J. Algebra \textbf{465} (2016), 259--286.

\bibitem[Bo02]{Bourbaki} N. Bourbaki. Lie groups and Lie algebras. Chapters 4--6. Springer, 2002.

\bibitem[BL]{BilleyLakshmibai} S. Billey, V. Lakshmibai. Singular loci of Schubert varieties. Progr. in Math. \textbf{182}, Birkh\"auser, 2000.



\bibitem[EP13]{EliseevPanov} D.Yu. Eliseev, A.N. Panov. Tangent cones to Schubert varieties for $A_n$ of lower rank.\break J. Math. Sci. \textbf{188} (2013), no. 5, 596--600.

\bibitem[EI13]{EliseevIgnatyev} D.Yu. Eliseev, M.V. Ignatyev. Kostant--Kumar polynomials and tangent cones to Schubert varieties for involutions in $A_n$, $F_4$ and $G_2$. J. Math. Sci. \textbf{199} (2014), no. 3, 289--301.

\bibitem[Hu175]{Humphreys} J. Humphreys. Linear algebraic groups. Springer, 1975.

\bibitem[Hu292]{Humpreys2} J. Humphreys. Reflection groups and Coxeter groups. Cambridge University Press, Cambridge, 1992.




\bibitem[IS16]{IgnatyevShevchenko1} M.V. Ignatyev, A.A. Shevchenko. On tangent cones to Schubert varieties in type $D_n$.\break St. Petersburg Math.~J. \textbf{27} (2016), no. 4, 609--623.



\bibitem[IS20]{IgnatyevShevchenko20} M.V. Ignatyev, A.A. Shevchenko. On tangent cones to Schubert varieties in type $E$. Comm. in Math. \textbf{28} (2020), 179--197.




\bibitem[Ki62]{Kirillov62} A.A. Kirillov. Unitary representations of nilpotent Lie groups. Russian Math. Surveys \textbf{17}~(1962), 53--110.

\bibitem[Ki03]{Kirillov03} A.A. Kirillov. Two more variations on the triangular theme. In: C. Duval, L. Guieu, V. Ovsienko, eds. The orbit method in geometry and physics. Progr. in Math. \textbf{213}, Birkh\"auser Boston, 2003, 243--258.

\bibitem[Ki04]{Kirillov04} A.A. Kirillov. Lectures on the orbit method. Grad. Studies in Math. \textbf{64}, AMS, 2004.







\bibitem[Ku96]{Kumar} S. Kumar. The nil-Hecke ring and singularities of Schubert varieties. Invent. Math. \textbf{123} (1996), 471--506.

\bibitem[Ku02]{Kumar02} S. Kumar. Kac--Moody groups, their flag varieties and representation theory. Progr. in Math. \textbf{204}, Birkh$\ddot{\mathrm{a}}$user, Boston, 2002.

\bibitem[Lu83]{Lusztig83} G. Lusztig. Some examples of square integrable representations of semisimple $p$-adic groups. Trans. Amer. Math. Soc. \textbf{277} (1983), 623--653.

\bibitem[Sh86]{Shi86} Y. Shi. The Kazhdan--Lusztig cells in certain affine Weyl groups. Lecture Notes in Math. \textbf{1179}. Springer--Verlag, Berlin, 1986





\end{thebibliography}
\end{document}